\documentclass[review]{elsarticle}
\usepackage{amsfonts,amssymb}
\usepackage{amsmath}
\usepackage{ntheorem}
\newtheorem{Theorem}{Theorem}
\newtheorem{Lemma}{Lemma}
\newtheorem*{Proof}{Proof}
\usepackage{cases}
\usepackage{bm}

\usepackage{threeparttable}
\biboptions{sort&compress}

\journal{Journal of \LaTeX\ Templates}

\begin{document}

\begin{frontmatter}

\title{Degree reduction of disk rational B\'{e}zier curves}

\author{Mao Shi\fnref{myfootnote}}
\address{School of Mathematics and Information Science of Shaanxi Normal University, Xi'an 710062, China}


\cortext[mycorrespondingauthor]{Corresponding author}
\ead{shimao@snnu.edu.cn}


\begin{abstract}
How to quickly and stably realize the degree reduction of the rational B\'{e}zier curve is an open problem in CAGD. Based on the weighted least squares method and weighted sum method of multi-objective optimization, this paper transforms the degree reduction problem of the rational B\'{e}zier curve into a convex optimization problem and then uses quadratic programming to solve it. Prove that the solution is the minimum. Numerical experiments show that the method is fast and stable.
\end{abstract}


\end{frontmatter}

%\linenumbers

\section{Introduction}
For the stability and robustness of numerical operation, the interval algorithm was brought into the geometric design system.
  In 1992, Sederberg and Farouki \cite{aSederberg} formally defined the interval B\'{e}zier curve.
  It can transfer a complete description of the approximate error along with the curve to other systems for the application.
  Thereafter, a series of practical algorithms such as curve/curve or surface/surface intersection, solid modeling, visualization stability, degree reduction, interval curve boundary problems, etc. were studied by \citep{Hu96a, Hu96b, Hu96c, flChen, Lin2002}.  However,  interval curves have some drawbacks. That is, the interval generally enlarges rapidly in a computational process and the rectangular intervals are not rotationally symmetric \cite{flChen}. However, Disk B\'{e}zier curves given by  Lin and Rokne \cite{Lin98} can correct these shortcomings. On the other hand, because  B\'{e}zier curves can't represent conic precisely, Hu et al. \cite{Hu96b}\cite{Hu96c} studied interval non-uniform rational B-splines (INURBS) curves and surfaces. In 2011, Based on parallel projection Shi \cite{Shi2011} defined a disk rational B\'{e}zier curve, which differs from the classic disk rational B\'{e}zier in that its error radii are B\'{e}zier polynomial functions.

Degree reduction of parametric curves is another important problem in the geometric modeling system.  It is useful in the data communication between the design system and data compression  \cite{Chen2004}. Compared with research on the degree reduction of the B\'{e}zier curve, the research on the degree reduction of the  rational B\'{e}zier curve is less\cite{flChen} \cite{baba2015}.  Farin \cite{Farin1983} described a degree reduction method for rational B\'{e}zier curves for interactive interpolation and approximation. Sederberg and Chang \cite{Sederberg1993}, Chen \cite{Chen1994} and Sun et al.\cite{Sun98} achieved  reduction through the approximate common divisor method. Qin and Guan \cite{Chen1993} transformed the problem of approximating multi-degree reduction of rational curves and surfaces into quadratic programming. Cai and Wang \cite{Cai} studied degree reduction of rational B\'{e}zier curves with the steepest descent method. Applying multi-objective optimization techniques, Shi \cite{Shi2015} realized multi-degree reduction rational B\'{e}zier curves.

Each of the above methods has its pros and cons. Some lack robustness in floating-point environments \cite{Sederberg1993} \cite{Chen1994} \cite{Sun98}; some have slower calculation speeds \cite{Cai}, etc. Therefore, finding a simple and stable method to obtain the optimal degree reduction approximation of rational curves is still an open problem. In this paper, multi-objective optimization, weighted least squares, and quadratic programming are serve to dealing with degree reduction of disk rational B\'{e}zier curves. Compared with the previous methods, this method is simple and stable.

The paper has the following structure: In section 2, we review the definition and properties of disk rational B\'{e}zier curve in  \cite{Shi2011}. In section 3, we propose an efficient algorithm to the problem of degree reduction of rational disk B\'{e}zier curves. In section 4, bounding errors for degree reduction are analyzed, and some examples are provided.

\section{Disk rational B\'{e}zier curves}
\subsection{Disk rational arithmetic}
Let $\mathbb{R}$  and $\mathbb{R}^{+}$  denote  the set of all real numbers and all nonnegative real numbers, respectively. A disk in the plane is defined as
\begin{equation}\label{eq:01}
  (\textbf{q})=(x_{0},y_{0})_{r}=\{\textbf{x}\in \mathbb{R}^{2}\ |\ \|\textbf{x}-\textbf{q}\|\leq r, r\in \mathbb{R}^{+}\ \},
\end{equation}
where  $\textbf{q}$ is  the centric point and $r$ is the radius.

For $n$ disks $(\textbf{q}_{i})=(x_{i},y_{i})_{r_{i}}$, we have the following linear combination \cite{flChen} \cite{Lin98}
\begin{eqnarray}\label{eq:02}
\sum_{i=0}^{n}k_{i}(\textbf{q}_{i})=\sum_{i=0}^{n}(k_{i}\textbf{q}_{i})=\left(\sum_{i=0}^{n}k_{i}x_{i},\sum_{i=0}^{n}k_{i}y_{i}\right)_{\sum\limits_{i=0}^{n}|k_{i}|r_{i}}, \ k \in \mathbb{R}.
\end{eqnarray}
Similarly, the disk in three-dimensional space can be  described as
\begin{eqnarray}\label{eq:03}
\begin{split}
(\textbf{P}^{\omega})&=(X_{0}, Y_{0}, \omega)_{r}=(\omega x_{0}, \omega y_{0}, \omega)_{ r}\\
&=\{ \textbf{x}^{\omega}=(\omega x,\omega y,\omega)\in \mathbb{R}^{3}\ |\  \|\textbf{\textit{x}}^{\omega}-\textbf{P}^{\omega} \|\leq r \}.
\end{split}
\end{eqnarray}
Applying the oblique projection $I(\cdot)$ to the disk $(\textbf{P}^{\omega})$,  we  obtain a  rational disk  in the plane $\omega=1$
\begin{eqnarray}\label{eq:04}
(\textbf{q})=I((\textbf{P}^{\omega}))=I(X_{0}, Y_{0}, \omega, r)=\left( \frac{X_{0}}{\omega}, \frac{Y_{0}}{\omega} \right )_{r}=(x_{0}, y_{0})_{r}.
\end{eqnarray}

\subsection{Disk rational B\'{e}zier curves}
On the basis of equations  \eqref{eq:02} and  \eqref{eq:04}, a disk rational curve can be defined, whose radius function is a positive polynomial function, and its properties  are the same as the disk B\'{e}zier one \cite{Shi2011}.
\\
\textbf{Definition 1.}\ \  A disk rational B\'{e}zier curve of degree $n$ with control disk points $({\textbf{p}}_{i})=(x_{i},y_{i})_{r_{i}}$  and corresponding weights $\omega_{i}\in \mathbb{R}^{+}, i=0,...,n,$ can be written as
\begin{eqnarray}\label{eq:05}
\left( \mathbf{p} \right) \left( t \right) =\left( \textbf{p}\left( t \right) ;r\left( t \right) \right) =\left( \frac{\sum\limits_{i=0}^n{\textbf{p}_i}\omega _iB_{i}^{n}\left( t \right)}{\sum\limits_{i=0}^n{\omega _i}B_{i}^{n}\left( t \right)};\sum_{i=0}^n{r_i}B_{i}^{n}\left( t \right) \right), 0\leq t\leq 1,
\end{eqnarray}
where $B_{i}^{n}(t)={n \choose i}t^{i}(1-t)^{n-i}, i=0,...,n,$ are Bernstein polynomials,  ${\textbf{p}}(t)$ and $r(t)$ are  called the center curve and the radius of the disk rational B\'{e}zier curve $(\textbf{{p}})(t)$, respectively.

Given another disk rational B\'{e}zier curves  $(\check{{\textbf{p}}})(t)$ of degree $m$ $(m<n)$
 \begin{equation}\label{eq:06}
  \left( {\check{\mathbf{p}}} \right) \left( t \right) =\left( {\check{\mathbf{p}}}\left( t \right) ;\check{r}\left( t \right) \right) =\left( \frac{\sum\limits_{i=0}^m{{\check{\mathbf{p}}}_i}\check{\omega}_iB_{i}^{m}\left( t \right)}{\sum\limits_{i=0}^m{\check{\omega}_i}B_{i}^{m}\left( t \right)};\sum_{i=0}^m{\check{r}_i}B_{i}^{m}\left( t \right) \right),
 \end{equation}
if the following equations  are satisfied
\begin{eqnarray}\label{eq:07}
\begin{array}{cc}
 \sum \limits_{j=max(0,i-n)}^{min(m,i) }{{m \choose j}{n \choose i-j}}\check{\omega}_{j}\omega_{i-j}\left({\textbf{p}}_{i-j}-\check{{\textbf{p}}}_{j}\right)= \mathbf{0}, \  i=0,1,...,n+m ,
 \end{array}
 \end{eqnarray}
 \begin{eqnarray}\label{eq:08}
r_{i}=\sum_{j=max(0,i-n+m)}^{min(m,i)}\frac{{m \choose j}{n-m \choose i-j}\check{r}_{j}}{{n \choose i}}, \ \  i=0,1,\ldots,n,
 \end{eqnarray}
 we say that  a degree $n$ disk rational B\'{e}zier curve $({\textbf{{p}}})(t)$ can be  represented exactly by a degree $m$ disk rational B\'{e}zier curve $(\check{\textbf{{p}}})(t)$.

 However, in many cases, we must use approximation methods to achieve the degree reduction of the disk rational B\'{e}zier curve.

 \section{Degree reduction of disk rational B\'{e}zier curves}

 \subsection{Description of approximation problem}

 The problem of  degree reduction of  disk Rational B\'{e}ziers curve can be stated as follows:
\\
\textbf{Problem 1.}  \ \ Given a degree $n$ disk rational B\'{e}zier curve $(\textbf{{p}})(t)$, find a degree $m<n$ disk rational B\'{e}zier curve $(\check{\textbf{{p}}})(t)$ such that $(\check{\textbf{{p}}})(t)$  is the closure of $(\textbf{{p}})(t)$.

By the weighted sum method of multi-objective optimization \cite{collette}, the above problem can be summarized as
\begin{subnumcases}
{}
(A)\ \min \ \  \frac{1}{2}\left( F_x+F_y \right)  \label{eq:09a}\\
   \ \ \ \  \ s.t. \   \  \check{\omega_{i}}>0, i=0,...,m \nonumber\\
(B) \min\limits_{\check{r}(t)\geq r(t)+dist(\textbf{{p}}(t),{\check{\textbf{{p}}}}(t))}\|\check{r}(t)-r(t)\|, \label{eq:09b}
 \end{subnumcases}
 where  $dist(\textbf{{p}}(t),\check{\textbf{{p}}}(t)),t\in [0, 1]$ is the Hausdorff distance between the curve $\check{{\textbf{p}}}(t)$ and the curve $\textbf{{p}}(t)$,
 $\rho(t)>0$ is a weight function and
\begin{equation}\label{eq:10}
  \left( F_x,\ F_y \right) =\int_0^1{\rho \left( t \right) \left( \textbf{p}\left( t \right) -\check{\textbf{{p}}}\left( t \right) \right)^2 dt}.
\end{equation}
That is, the degree reduction of the disk rational b\'{e}zier curve  is composed of two parts, one is the degree reduction of the central curve $\check{{\textbf{p}}}(t)$, and the other for  the error radius curve $r(t)$.
\subsection{Degree reduction approximation of center curve}

In order to ensure that all weights  are positive,  we give some basic theorems firstly and then use quadratic programming  to achieve the degree reduction of the central curve.
\begin{Lemma}{\rm\cite{Cai}}
  The rational B\'{e}zier curves $\check{\textbf{p}}(t)$ and $\textbf{p}(t)$ satisfy $C^{(u, v)}$-continuity if and only if the following equations are true:
  \begin{equation}\label{eq:11}
    \sum_{i=0}^q{{n \choose {q-i}}{m\choose i} \omega _{q-i} \check{\omega}_i\check{\textbf{{p}}}_i}=\sum_{i=0}^q{{m \choose {q-i}}{n\choose i} \check{\omega}_{q-i} \omega _i\textbf{p}_i},\ q=0,\cdots ,u;
  \end{equation}
  and
\begin{align}\label{eq:12}
  &\sum_{i=0}^l{{n \choose {n-i}}{m\choose m-l+i} \omega _{n-i} \check{\omega}_{m-l+i}\check{\textbf{{p}}}_{m-l+i}}
  =\sum_{i=0}^l{{m \choose {m-i}}{n\choose n-l+i} \check{\omega}_{m-i} \omega _{n-l+i}\textbf{p}_{n-l+i}},\nonumber\\
  &l=0,\cdots ,v.
\end{align}
\end{Lemma}
Rewriting Lemma 1, we have
\begin{Theorem}
Any component of the vector
$
\mathbf{\check{P}}^{\mathbf{I}}=\left[ \check{\omega}_0\mathbf{\check{p}}_0,\cdots ,\check{\omega}_u\mathbf{\check{p}}_u \right] ^{\rm T}
$
can be  expressed as a linear combination of components in the vector
$
\bm{\check{\omega} }^{\mathbf{I}}=\left[ \check{\omega}_0,\cdots ,\check{\omega}_u \right] ^{\rm T}.
$
It also holds for
$
\mathbf{\check{P}}^{\mathbf{III}}=\left[ \check{\omega}_{m-v}\mathbf{\check{p}}_{m-v},\cdots ,\check{\omega}_m\mathbf{\check{p}}_m \right] ^{\rm T}
$
 and $
\bm{\omega }^{\mathbf{III}}=\left[ \check{\omega}_{m-v},\cdots ,\check{\omega}_m \right] ^{\rm T}.
$
That is
\begin{equation}\label{eq:13}
  \mathbf{\check{P}}^{\mathbf{I}}=\mathbf{S}_{1}^{-1}\mathbf{T}_1\bm{\check{\omega} }^{\mathbf{I}}
\end{equation}
and
\begin{equation} \label{eq:14}
  \mathbf{\check{P}}^{\mathbf{III}}=\mathbf{S}_{3}^{-1}\mathbf{T}_3\bm{\check{\omega} }^{\mathbf{III}},
\end{equation}
where
\\
$
\mathbf{S}_1\!=\!\left[s^1_{ij}\right]\!=\!\left[ {n\choose i-j}{m\choose j} \omega _{i-j} \right],
\mathbf{T}_1\!=\!\left[t^1_{ij}\right]\!=\!\left[{n\choose i-j}{m\choose j}  \omega _{i-j}\mathbf{p}_{i-j} \right],
(0\leq j \leq i \leq u),$
$
\mathbf{S}_3\!=\!\left[s^3_{kl}\right]\!=\!\left[ {n\choose n+k-l}{m\choose m-v+l}\omega _{n+k-l} \right],
\mathbf{T}_3\!=\!\left[t^3_{kl}\right]\!=\!\left[{n\choose n+k-l}{m\choose m-v+l} \omega _{n+k-l}\mathbf{p}_{n+k-l} \right],
(0\leq k\leq l\leq v).
$
\end{Theorem}
\begin{Proof}
  Equation \eqref{eq:11} can be rewritten in matrix form
   \begin{align*}
    &\left[ \begin{smallmatrix}
	{n\choose 0}\omega _0&		0&		\cdots&		0\\
	{n\choose 1} \omega _1&		{n\choose 0} \omega _0&		\cdots&		0\\
	\vdots&		\vdots&		\ddots&		\vdots\\
	{n\choose u} \omega _u&		{n\choose u-1} \omega _{u-1}&		\cdots&		{n\choose 0} \omega _0\\
\end{smallmatrix} \right] \left[\begin{smallmatrix}
	{m\choose 0}&		0&		\cdots&		0\\
	0&		{m\choose 1}&		\cdots&		0\\
	\vdots&		\vdots&		\ddots&		\vdots\\
	0&		0&		\cdots&		{m\choose u}\\
\end{smallmatrix} \right] \left[\begin{smallmatrix}
	\check{\omega}_0\mathbf{\check{p}}_0\\
	\check{\omega}_1\mathbf{\check{p}}_1\\
	\vdots\\
	\check{\omega}_u\mathbf{\check{p}}_u\\
\end{smallmatrix} \right] \\
&=
\left[ \begin{smallmatrix}
	{n\choose 0}\omega _0\mathbf{p}_0&		0&		\cdots&		0\\
	{n\choose 1}\omega _1\mathbf{p}_1&	{n\choose 0}	 \omega _0\mathbf{p}_0&		\cdots&		0\\
	\vdots&		\vdots&		\ddots&		\vdots\\
	{n\choose u} \omega _u\mathbf{p}_u&		{n\choose u-1} \omega _{u-1}\mathbf{p}_{u-1}&		\cdots&	{n\choose 0}	 \omega _0\mathbf{p}_0\\
\end{smallmatrix} \right] \left[ \begin{smallmatrix}
	{m\choose 0}&		0&		\cdots&		0\\
	0&	{m\choose 1}	&		\cdots&		0\\
	\vdots&		\vdots&		\ddots&		\vdots\\
	0&		0&		\cdots&		{m\choose u}\\
\end{smallmatrix} \right] \left[ \begin{smallmatrix}
	\check{\omega}_0\\
	\check{\omega}_1\\
	\vdots\\
	\check{\omega}_u\\
\end{smallmatrix} \right]
   \end{align*}
   Regrouping the above matrix equation, it yields
\begin{align*}
&\left[ \begin{smallmatrix}
	\check{\omega}_0\mathbf{\check{p}}_0\\
	\check{\omega}_1\mathbf{\check{p}}_1\\
	\vdots\\
	\check{\omega}_u\mathbf{\check{p}}_u\\
\end{smallmatrix} \right] =\left[ \begin{smallmatrix}
	{n \choose 0}{m\choose 0}\omega _0&		0&		\cdots&		0\\
	{n \choose 1}{m\choose 0} \omega _1&	{n \choose 0}{m\choose 1} \omega _0&		\cdots&		0\\
	\vdots&		\vdots&		\ddots&		\vdots\\
	{n \choose u}{m\choose 0} \omega _u&		{n \choose u-1}{m\choose 1}\omega _{u-1}&		\cdots&	{n \choose 0}{m\choose u}	 \omega _0\\
\end{smallmatrix} \right] ^{-1}\\
&\times\left[\begin{smallmatrix}
	{m\choose 0}{n\choose 0} \omega _0\mathbf{p}_0&		0&		\cdots&		0\\
	{n\choose 1}{m\choose 0} \omega _1\mathbf{p}_1&	 {n\choose 0}{m\choose 1}	 \omega _0\mathbf{p}_0&		\cdots&		0\\
	\vdots&		\vdots&		\ddots&		\vdots\\
	{n\choose u}{m\choose 0} \omega _u\mathbf{p}_u&	{n\choose u-1}{m\choose 1}	\omega _{u-1}\mathbf{p}_{u-1}&		\cdots&	{n\choose 0}{m\choose u}\omega _0\mathbf{p}_0\\
\end{smallmatrix}\right]\left[ \begin{smallmatrix}
	\check{\omega}_0\\
	\check{\omega}_1\\
	\vdots\\
	\check{\omega}_u\\
\end{smallmatrix} \right],
\end{align*}
which establishes  equation \eqref{eq:11}. Similarly, equation \eqref{eq:14} also holds by \eqref{eq:12}.
\end{Proof}
Letting
\begin{equation} \label{eq:15}
  \rho \left( t \right) =\left( \sum\limits_{i=0}^m{\check{\omega}_iB_{i}^{m}\left( t \right)}\sum\limits_{i=0}^n{\omega _iB_{i}^{n}\left( t \right)} \right) ^2
\end{equation}
 and  substituting  $\mathbf{p}\left( t \right) $ and $\check{\bf{p}}\left( t \right)$ into  equation \eqref{eq:10}, it yields
  \begin{align} \label{eq:16}
   &\left[ F_x\left( t \right) ,F_y\left( t \right) \right] \nonumber\\
   &=\int_0^1{\left( \sum_{i=0}^m{\check{\omega}_iB_{i}^{m}\left( t \right)}\sum_{i=0}^n{\omega _i\mathbf{p}_iB_{i}^{n}\left( t \right)}-\sum_{i=0}^m{\check{\omega}_i\mathbf{\check{p}}_iB_{i}^{m}\left( t \right) \sum_{i=0}^n{\omega _iB_{i}^{n}\left( t \right)}} \right) ^2dt}.
  \end{align}
 When the above equation reaches its minimum, we must have
  \begin{align*}
    &\frac{\partial \left( F_x\left( t \right) ,F_y\left( t \right) \right)}{\partial \mathbf{\check{p}}_k}  \ \ \ \  (k=u+1, \ldots, m-v-1) \nonumber\\
    &\!=\!2\int_0^1{\left(\! \sum_{i=0}^n{\omega _i\mathbf{p}_iB_{i}^{n}\left( t \right)}\!\sum_{i=0}^m{\check{\omega}_iB_{i}^{m}\left( t \right)}\!-\!\sum_{i=0}^n{\omega _iB_{i}^{n}\left( t \right)}\sum_{i=0}^m{\check{\omega}_i\mathbf{\check{p}}_iB_{i}^{m}\left( t \right)} \right) \sum_{i=0}^n{\omega _iB_{i}^{n}\left( t \right)}}B_{k}^{m}\left( t \right) dt \nonumber\\
    &=\mathbf{0}
  \end{align*}
 Unfolding the above equation and writing it in matrix form, we deduce
 \begin{equation}\label{eq:17}
\mathbf{S}_2\mathbf{\check{P}}^{\mathbf{II}}=\mathbf{D}\bm{\check{\omega} }-\mathbf{S}_{2}^{\mathbf{L}}\mathbf{\check{P}}^{\mathbf{I}}-\mathbf{S}_{2}^{\mathbf{R}}\mathbf{\check{P}}^{\mathbf{III}},
\end{equation}
where $
\mathbf{\check{P}}^{\mathbf{II}}=\left[ \check{\omega}_{u+1}\mathbf{\check{p}}_{u+1},\cdots ,\check{\omega}_{m-v-1}\mathbf{\check{p}}_{m-v-1} \right] ^{\rm T},
$
$
\bm{\check{\omega} }=\left[ \check{\omega}_{0},\cdots ,\check{\omega}_{m} \right] ^{\rm T},
$
$$
\mathbf{S}_2=\left[ S_{ij}^{2} \right] =\left[ {\sum\limits_{h = 0}^{2n} {\frac{{{m\choose {u+i}}{2n\choose h}{m\choose {u+j}}A_h }}{{\left( {2m + 2n + 1} \right){{2m+2n}\choose {2u+i+j+h}}}}} } \right],\ (i,j=1,\ldots,m-u-v-1),
$$
$$
\mathbf{D}=\left[ D_{ij} \right] =\left[ {\sum\limits_{h = 0}^{2n} {\frac{{{m\choose {u+i}}{2n\choose h}{m\choose {j}}\mathbf{C}_h }}{{\left( {2m + 2n + 1} \right){{2m+2n}\choose {u+i+j+h}}}}} } \right],\ (i=1,\ldots,m-u-v-1;j=0,\ldots,m),
$$
$$
\mathbf{S}_{2}^{\mathbf{L}}=\left[ S_{ij}^{L} \right] =\left[ {\sum\limits_{h = 0}^{2n} {\frac{{{m\choose {u+i}}{2n\choose h}{m\choose {j}}A_h }}{{\left( {2m + 2n + 1} \right){{2m+2n}\choose {u+i+j+h}}}}} } \right], \
(i=1,\ldots,m-u-v-1;j=0,\ldots,u),
$$
$$
\mathbf{S}_{2}^{\mathbf{R}}=\left[ S_{ij}^{R} \right] =\left[ {\sum\limits_{h = 0}^{2n} {\frac{{{m\choose {u+i}}{2n\choose h}{m\choose {m-v+j}}A_h }}{{\left( {2m + 2n + 1} \right){{2m+2n}\choose {m-v+u+i+j+h}}}}} } \right], \
(i=1,\ldots,m-u-v-1;j=0,\ldots,v),
$$
$A_i=\sum\limits_{j=\max \left( 0,i-n \right)}^{\min \left( n,i \right)}{\frac{{n\choose j}{n\choose i-j} }{{2n\choose i}}\omega _j\omega _{i-j}}$  and $ \mathbf{C}_i=\sum\limits_{j=\max \left( 0,i-n \right)}^{\min \left( n,i \right)}{\frac{{n\choose j}{n\choose i-j} }{{2n\choose i}}\omega _j\mathbf{p}_j\omega _{i-j}}.$
\\
\\
Substituting equations \eqref{eq:13} and \eqref{eq:14} into equation \eqref{eq:17}, we obtain
\begin{Theorem}
Under the action of the  weighted least squares  and $\rho \left( t \right)$, any component of the vector
$
\mathbf{\check{P}}^{\mathbf{II}}$
 is a linear combinations of components in the vector $
\bm{\check{\omega} },$
 and
 \begin{align}\label{eq:18}
   \mathbf{\check{P}}^{\mathbf{II}} & =\mathbf{S}_{2}^{-1}\left( \mathbf{D}\bm{\check{\omega}}-\mathbf{S}_{2}^{\mathbf{L}}\mathbf{S}_{1}^{-1}\mathbf{T}_1\bm{\check{\omega}}^{\mathbf{I}}-\mathbf{S}_{2}^{\mathbf{R}}\mathbf{S}_{3}^{-1}\mathbf{T}_3\bm{\check{\omega}}^{\mathbf{III}} \right) \nonumber\\
    &=\mathbf{S}_{2}^{-1}\left( \mathbf{D}-\left[ \begin{matrix}
	\mathbf{S}_{2}^{\mathbf{L}}\mathbf{S}_{1}^{-1}\mathbf{T}_1&		\mathbf{0}_{\left( m-u-v-1 \right) \times \left( m-u-v-1 \right)}&		\mathbf{S}_{2}^{\mathbf{R}}\mathbf{S}_{3}^{-1}\mathbf{T}_3\\
\end{matrix} \right] \right) \bm{\check{\omega}}.
 \end{align}
\end{Theorem}
Finally, combining \textbf{Theorems} 1 and 2, it yields
\begin{Theorem}
 $\mathbf{\check{P}}$ and $\bm{\check{\omega}}$ have the following relation
\begin{equation}\label{eq:19}
\mathbf{\check{P}}=\mathbf{Q}\bm{\check{\omega}},
\end{equation}
where $\mathbf{\check{P}}=\left[\check{\omega}_0\mathbf{\check{p}}_0 ,		\cdots ,		\check{\omega}_m\mathbf{\check{p}}_m
\right] ^{\rm T}$ and
$$
\mathbf{Q}=\left[ \begin{array}{cc}
	\mathbf{S}_{1}^{-1}\mathbf{T}_1&		\mathbf{0}_{\left( u+1 \right) \times \left( m-u \right)}\\
	\multicolumn{2}{c}{\mathbf{S}_{2}^{-1}\left( \mathbf{D}-\left[ \begin{matrix}
	\mathbf{S}_{2}^{\mathbf{L}}\mathbf{S}_{1}^{-1}\mathbf{T}_1&		\mathbf{0}_{\left( m-u-v-1 \right) \times \left( m-u-v-1 \right)}&		\mathbf{S}_{2}^{\mathbf{R}}\mathbf{S}_{3}^{-1}\mathbf{T}_3\\
\end{matrix} \right] \right) }		\\
	\mathbf{0}_{\left( v+1 \right) \times \left( m-v \right)}&		\mathbf{S}_{3}^{-1}\mathbf{T}_3\\
\end{array} \right].
$$
\end{Theorem}
Using equation \eqref{eq:19},  equation \eqref{eq:16} can be expressed as

\begin{equation} \label{eq:20}
\left[F_x, F_y \right]=\bm{\check{\omega}}^{\rm T}\int_0^1{\left( \mathbf{P}^\text{T}\mathbf{K}-\bm{\omega}^{\rm T}\mathbf{ KQ} \right)^{\rm T}\left( \mathbf{P}^\text{T}\mathbf{K}-\bm{\omega}^{\rm T}\mathbf{ KQ} \right) dt}\bm{\check{\omega}},
\end{equation}
where $
\mathbf{P}=\left[
	\omega _0\boldsymbol{p}_0,		\cdots,		\omega _n\boldsymbol{p}_n
 \right] ^{\rm T},\ \bm{\omega }=\left[ 	\omega _0,		\cdots,		\omega _n
\right] ^{\rm T}
$ and $$
\mathbf{K}=\left[ \begin{matrix}
	B_{0}^{n}\left( t \right) B_{0}^{m}\left( t \right)&		\cdots&		B_{0}^{n}\left( t \right) B_{m}^{m}\left( t \right)\\
	\vdots&		\ddots&		\vdots\\
	B_{n}^{n}\left( t \right) B_{0}^{m}\left( t \right)&		\cdots&		B_{n}^{n}\left( t \right) B_{m}^{m}\left( t \right)\\
\end{matrix} \right].
$$
Obviously,  equation \eqref{eq:20} is positive, so we have
\begin{Theorem}
 The Hessian matrix of equation \eqref{eq:16}
 \begin{equation}\label{eq:21}
\left[H_x, H_y\right]=\int_0^1{\left( \mathbf{P}^{\rm T}\mathbf{K}-\bm{\omega}^{\rm T}\mathbf{ KQ} \right) ^{\rm T}\left( \mathbf{P}^{\rm T}\mathbf{K}-\bm{\omega}^{\rm T}\mathbf{ KQ} \right) dt}
\end{equation}
is positive definite, and the solution of equation \eqref{eq:09a} with respect to $\check{\omega}_i, \ (i=0,\ldots, m)$  is unique.
\end{Theorem}
Finally, $\bm{\check{\omega}}$ can be obtained by
\begin{align}\label{eq:22}
  \min \ \ \ & \frac{1}{2}\bm{\check{\omega}}^T\left( H_x+H_y \right) \bm{\check{\omega}} \\
   s.t. \ \ \ & \bm{\check{\omega}}>\bm{0},   \nonumber
\end{align}
and  $\left\{ \mathbf{\check{p}}_i \right\} _{i=0}^{m}$  by
$$
\mathbf{\check{p}}=\left[
	\mathbf{\check{p}}_0,		\cdots,		\mathbf{\check{p}}_m \right] ^\text{T}=\mathbf{M}^{-1}\mathbf{Q}\boldsymbol{\check{\omega}},
$$
where
$$
\mathbf{M}=diag \left[\omega_i, ..., \omega_m  \right].
$$
\section{The error function}
For  simplicity, we use the following function as the metric  degree reduction error function
\begin{equation}\label{eq:22}
 d=\int_0^1{\left( \mathbf{p}\left( t \right) -\mathbf{\check{p}}\left( t \right) \right) \cdot \left( \mathbf{p}\left( t \right) -\mathbf{\check{p}}\left( t \right) \right) dt},
\end{equation}
where multiplication between vectors is inner product.

\section{Degree reduction approximation of error radius curve}

Similar to the derivation of the degree reduction of error radius in \cite{Shi2015} , we have
\begin{eqnarray}\label{eq:32}
 \left \{
\begin{array}{cc}
\ min\  &\sum_{i=0}^{m}\sum_{j=0}^{m}\check{r}_{i}\check{r}_{j}H_{ij}-2\sum_{i=0}^{m}\sum_{j=0}^{n}\check{r}_{i}r_{j}S_{ij}
\\
s.t. \   &\hat{r}_{i}\geq r_{i}+d, \ \ \ i=0,1,\ldots, n,\\
 \ \ &\check{r}_{j}>0,  \  \   \  \   \  \   \   \  \  \        j=0,1,...,m,
 \end{array}
 \right .
\end{eqnarray}
where $H_{ij}=\left[\frac{{m \choose i}{m \choose j}}{(2m+1){2m \choose i+j}}\right]\ (i,j=0,\ldots,m)$, $S_{ij}=\left[\frac{{m \choose i}{n \choose j}}{(m+n+1){m+n \choose i+j}}\right] \ (i=0,\ldots,m;j=0, \ldots, n)$,
$\hat{r}_{i}$ are given by \eqref{eq:08} and $d$ is defined by \eqref{eq:22}.

\section{Numerical examples}
In this section, we give several examples to illustrate the effectiveness of our method. Except for the results of Cai and Wang\cite{Cai} used in Example 1, the results of other examples were obtained by Matlab2018b, which is somewhat different from the results given in \cite{Cai}.

{{\textbf{Example 1.} (Also Example 1 in \cite{Cai} and \cite{Sederberg1993}) Given a 4 degree rational B\'{e}zier curve with control points in homogenous coordinates:$
\left( \omega _i\boldsymbol{R}_i, \omega_i\right) =\left( 0,0,1 \right) ,\left( 8,8,4 \right) ,\left( 6,0,2 \right) ,\left( 4,-2,1 \right) ,\left( 4,0,1 \right)$, to find a 1-degree reduced rational B\'{e}zier curve to approximate the original curve.  Table 1 gives the values under the different error measures. The resulting curves are illustrated in Figure 1.
\begin{figure*}[htp]
   \centering
    %\begin{tabular}{cc}
    %\begin{minipage}[t]{3in}
    \includegraphics[width=4.in]{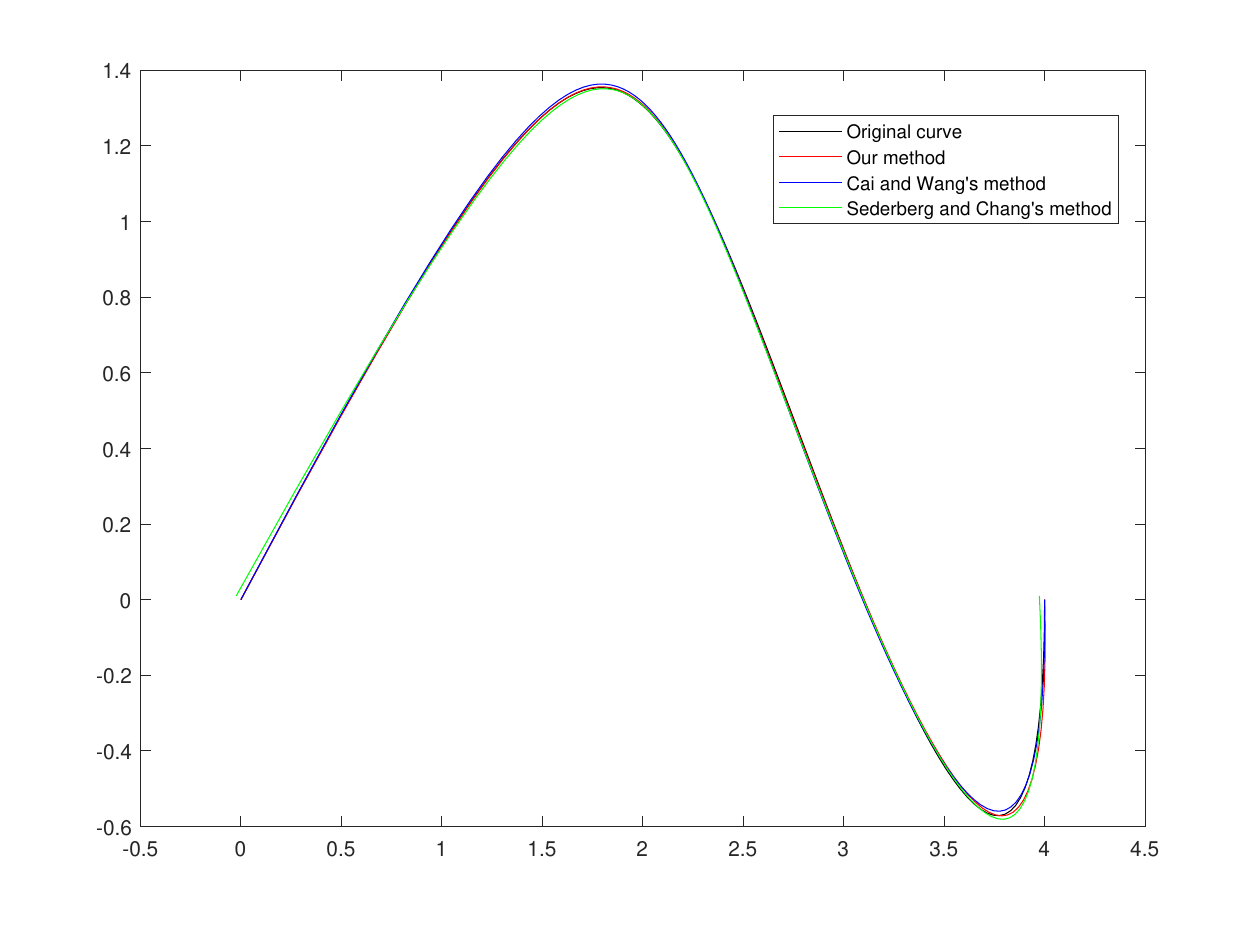}
    \caption{Comparison of the three degree reduction methods}\label{fig1}
    %\end{minipage}
    %\end{tabular}
\end{figure*}

\begin{table}[]
\caption{Error comparison of three methods.}
\begin{center}
\begin{tabular}{|l|l|}
\hline
Methods                      & Errors   \\ \hline
Our method                   & 0.007330 \\ \hline
Cai and Wang's method        & 0.008324 \\ \hline
Sederberg and Chang's method & 0.012064 \\ \hline
\end{tabular}
\end{center}
\end{table}

{{\textbf{Example 2.} (Also Example 2 in \cite{Cai}) Given a 5 degree rational B\'{e}zier curve with control points in homogenous coordinates:$
\left( \omega _i\boldsymbol{R}_i,\ \omega _i \right) =\left( 0,0,1 \right) ,\left( 4,20,2 \right)$, $\left( 24,48,4 \right) ,\left( 70,56,7 \right) ,\left( 14,2,2 \right) ,\left( 18,6,3 \right)$, to find a 1-degree reduced rational B\'{e}zier curve. Table 2 gives comparisons of approximation error and time. The resulting curves are illustrated in Figure 2.
\begin{figure*}[htp]
   \centering
    %\begin{tabular}{cc}
    %\begin{minipage}[t]{3in}
    \includegraphics[width=3.5in]{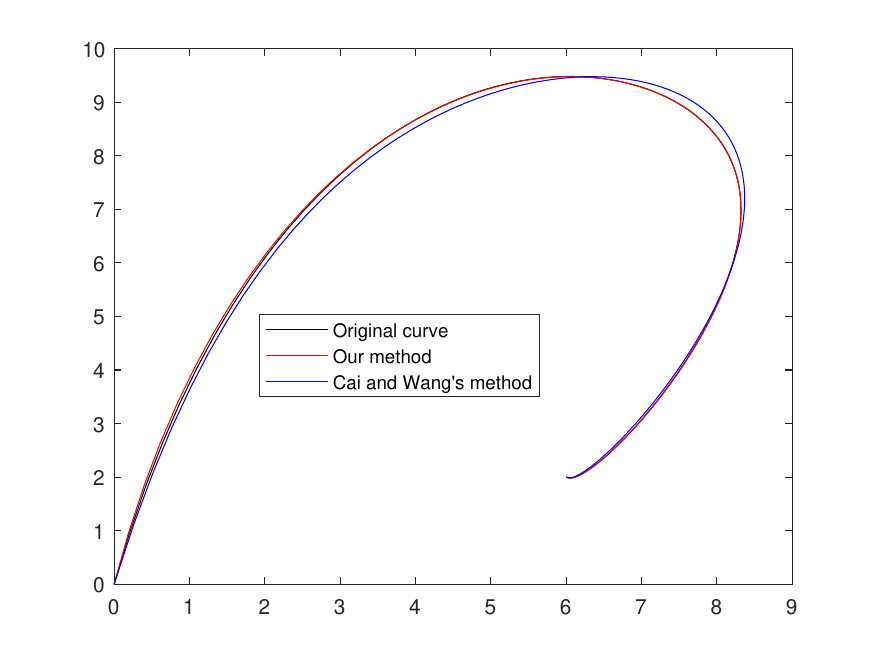}
    \caption{Comparison of the two methods }\label{fig2}
    %\end{minipage}
    %\end{tabular}
\end{figure*}

\begin{table}[]
\caption{Time and error comparisons of the two methods.}
\centering
\begin{threeparttable}
\begin{tabular}{|l|l|l|}
\hline
Methods             & Time (s)         & Error  \\ \hline
Our method            & 0.917970         & 0.0096 \\ \hline
Cai And Wang's method$^{*}$ & 14.455630 ($M$=15) & 0.1469 \\ \hline
\end{tabular}
\begin{tablenotes}
        \footnotesize
        \item[*] The results of CAI and Wang's method are given by our program. It is approximate to the error given in \cite{Cai}.  %此处加入注释*信息
        %此处加入注释**信息
      \end{tablenotes}
    \end{threeparttable}
\end{table}

{{\textbf{Example 3.} (Also Example 3 in \cite{Cai}) Given a 8 degree rational B\'{e}zier curve with control points in homogenous coordinates: $
\left( \omega _i\boldsymbol{R}_i,\ \omega _i \right) =\left( 0,0,1 \right) ,\left( 0,4,2 \right)$, $\left( 6,30,3 \right) ,\left( 36,54,9 \right) ,\left( 72,72,12 \right) ,\left( 220,320,20 \right) ,\left( 240,30,30 \right) ,\left( 36,4,4 \right) ,\left( 10,0,1 \right)$, to find a 3-degree reduced rational B\'{e}zier curve. Table 3 gives comparisons of approximation error and time. The resulting curves are illustrated in Figure 3.
\begin{figure*}[htp]
   \centering
    %\begin{tabular}{cc}
    %\begin{minipage}[t]{3in}
    \includegraphics[width=3.5in]{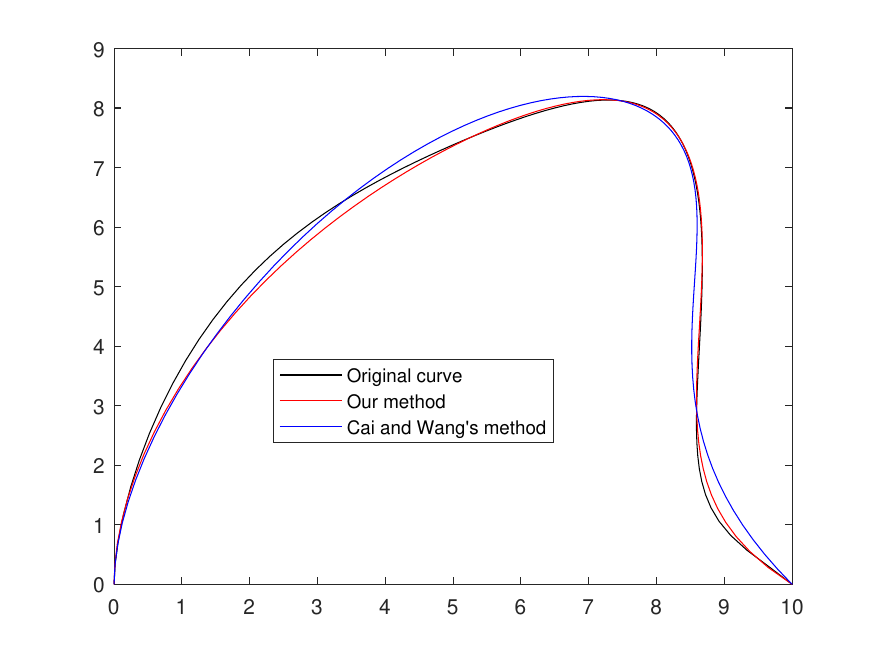}
    \caption{Comparison of the two methods }\label{fig3}
    %\end{minipage}
    %\end{tabular}
\end{figure*}

\begin{table}[]
\caption{Time and error comparisons of the two methods.}
\centering
\begin{threeparttable}
\begin{tabular}{|l|l|l|}
\hline
Methods             & Time (s)         & Error  \\ \hline
Our method            & 1.719172         & 0.1687\\ \hline
Cai And Wang's method$^{*}$ & 4.440010 ($M$=4) & 0.2358 \\ \hline
\end{tabular}
\begin{tablenotes}
        \footnotesize
        \item[*] The results of CAI and Wang's method are given by our program. It is approximate to the error given in \cite{Cai}.
      \end{tablenotes}
    \end{threeparttable}
\end{table}

{{\textbf{Example 4.} (Also Example 2 in \cite{Shi2015}) Given a 8 degree disk rational B\'{e}zier $(\textbf{p})(t)$  with control disks
$(6, 14.9)_1$, $(8.6, 25)_{0.4}$,$(20.3, 30)_1$, $(35, 31)_{1.5}$, $(40.2, 25)_2$, $(37.5, 11.5)_{1.8}$, $(47.2, 8.1)_{0.8}$, $(65.1, 11.2)_1$, $(71.5, 25)_{0.5}$ and associated weights
 1.88, 1.68, 1.63, 1.73, 1.79, 2.18, 1.24, 1.08, 1.9. The best 3-degree reduction curve satisfying $C^{(1,1)}$-continuity with the given curve has control
disks $(6.0000, 14.9000 )_{10.4812}$, $(11.2981, 35.4812)_{0.8812}$, $(56.8670, 34.1932)_{26.7312}$, $(31.7728, 13.5463)_{18.6645}$, $(48.8560, -23.8262)_{13.4812}$, $( 71.5000, 25.0000)_{5.4812}$, and
associated weights $0.0023,    0.0016,    0.0013,    0.0020,    0.0005,    0.0019$. The resulting curve is illustrated in Figure 4.
\begin{figure*}[htp]
   \centering
    %\begin{tabular}{cc}
    %\begin{minipage}[t]{3in}
    \includegraphics[width=3.5in]{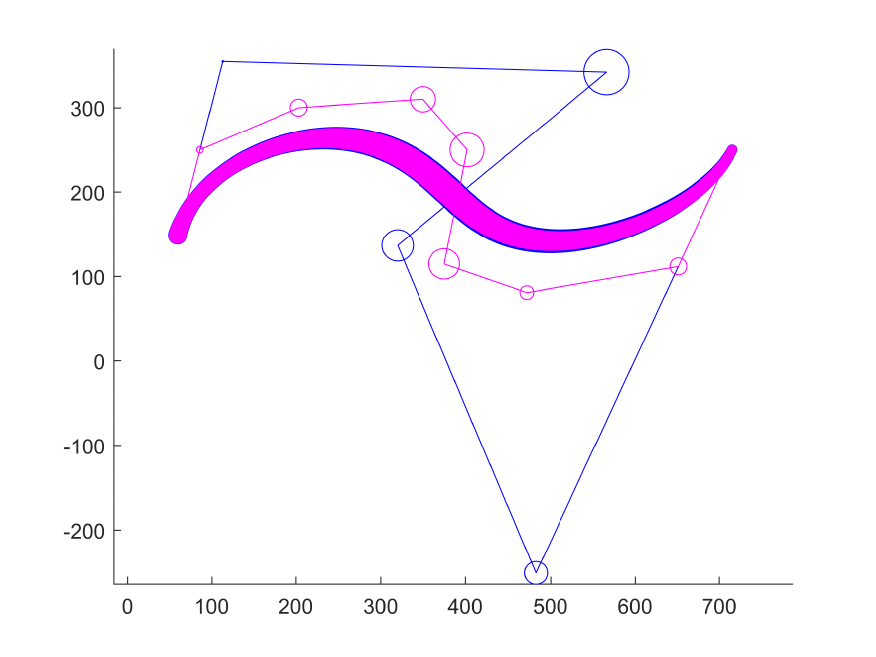}
    \caption{The resulting curve}\label{fig1}
    %\end{minipage}
    %\end{tabular}
\end{figure*}

\section*{References}

\bibliography{test}%{mybibfile}

\end{document}